\documentclass[10pt,letterpaper]{amsart}

\usepackage{amsfonts,amssymb,amsmath,amsthm}
\numberwithin{equation}{section}

\def\tfr#1/#2{{#1}/{#2}}			
\def\fr#1/#2{\frac{#1}{#2}}			

\newcommand\suny{\sum_{n=1}^\infty}
\newcommand\sumy{\sum_{m=1}^\infty}
\newcommand\suly{\sum_{l=1}^\infty}
\newcommand\sumn{\sum_{m=1, m\mid n}^n}

\author{Jacques G\'elinas}

\address{Ottawa, Canada}

\email{jacquesg00@hotmail.com}

\thanks{This work was done in May 2013 and July 2014 while the author was a retired mathematician}

\begin{document}

\title[Proof of RH by Lee]{Note on a proposed proof of the Riemann Hypothesis by Jin Gyu Lee}

\begin{abstract}

	This is a reformulation and refutation of a proposed proof of the Riemann hypothesis
	published in electronic form on the Internet in 2013 and updated in 2014.
	Proceeding by contradiction, the author wants to show that if $\zeta(s)=0$ where 
	$\tfr 1/2<\Re s<1$, then $\zeta(2s)=0$, which is known to be impossible. We show 
	that both versions of the proposed proof are incomplete.
\end{abstract}

\maketitle


\section{On the right of the critical strip}
Consider the $\zeta(s)$ function of Riemann, represented by a Dirichlet series, an Euler product,
or via the Dirichlet alternating zeta integral function $\eta(s)$ as~:
\begin{align}
	\zeta(s) &= \suny \, \fr 1/{n^s},	&(\Re s >1),
\\ \label{Eulerproduct}
		&=	\prod_{p \text{ prime}} \, \fr 1/{1-p^{-s}},	&(\Re s >1),
\\ \label{zeta:eta}
		&=	\fr 1/{1-2^{1-s}}\;\eta(s),	&(\Re s >0, s \ne 1+\fr {2k\pi i}/{\log 2}, k \text { integer}),
\\ \nonumber
		&=	\fr {\eta'(s)}/{\log 2},
				&(s = 1+\fr {2k\pi i}/{\log 2}, k \text{ nonzero integer}),
\\ \nonumber
	\eta(s) &= \suny \, \fr (-1)^{n+1}/{n^s},	&(\Re s >0).
\end{align}
In his second proof of the functional equation of the zeta function discovered by Euler \cite{Euler:1749},
Riemann \cite{Riemann:1859} has shown that $\xi(s):=\pi^{\tfr s/2}\Gamma(1+\tfr s/2)(s-1)\zeta(s)$ 
is an integral function left unchanged by the mapping $s \mapsto 1-s$. The non trivial
zeros of $\zeta(s)$ are thus located symmetrically with respect to the line $\Re s=\tfr 1/2$, inside the
critical strip $0<\Re s<1$ because the Euler product (\ref{Eulerproduct}) is not zero for $\Re s>1$, 
nor is the limit for $\Re s=1$ of an expression derived from it. Riemann has stated that all 
these non trivial zeros are very likely located on the critical line $\Re s=\tfr 1/2$ itself.

In order to prove this ``Riemann Hypothesis", it is sufficient to show that if $\eta(s)$ had one 
zero in the right hand side $\tfr 1/2 <\Re s<1$ of the critical strip,
then $\zeta(2s)$ would also vanish while $2s$ is outside the critical strip, contradicting
(\ref{Eulerproduct}). This is Mr Lee's claim and the main goal of \cite{Lee:2013,Lee:2014}~:
\begin{equation} \label{Leeclaim}
	\fr 1/2<\Re s <1 \;\&\; \suny \, \fr (-1)^{n+1}/{n^s}=0 
		\;\implies\;
			 \Re(2s)>1 \;\&\; \suny \fr 1/{n^{2s}} = 0.
\end{equation}

\section{Arithmetic functions}
The following two arithmetic functions are defined on the natural numbers in \cite{Lee:2013}~: first,
with the convention that $\Omega(1)=0$, \cite[definition 3.1]{Lee:2013}
\begin{equation} \label{defomeg}
	\Omega(m) := \sum_{k=0}^K r_k \quad\text{where, for distinct primes } p_k,\, m = \prod_{k=0}^K p_k^{r_k},
\end{equation}
and secondly \cite[definition 3.3]{Lee:2013}
\begin{equation} \label{defbeta}
	\beta(n) := \sumn (-1)^{\Omega(m)}(-1)^{l + 1}, \qquad(l=\fr n/m).
\end{equation}
According to \cite[theorem 3.7]{Lee:2013}~:
\begin{equation} \label{thmbeta}
	\beta(n)=
\begin{cases}
	1	&	\text{if $n$ is the square of a natural number},\\
	-2	&	\text{if $n$ is twice the square of a natural number},\\
	0	&	\text{otherwise}.
\end{cases}
\end{equation}

\section{A new expression for the Riemann zeta function}
If $\Re s>\tfr 1/2$, then $\Re (2s)>1$ and the Dirichlet series representing $\zeta(2s)$ converges
absolutely (and also uniformly on compact sets). We can certainly write, replacing the usual exponent
$1-2s$ for $\eta(2s)$ in (\ref{zeta:eta}) by the exponent $1-s$ as in \cite{Lee:2013}~:
\begin{align*}
	(1-2^{1-s}) & \suny \fr 1/{n^{2s}} \,=\,
			 \suny \, \left( \, \fr 1/{(n^2)^s} \,-\, \fr 2/{(2n^2)^s} \, \right),
\\
		&= \left( \fr 1/{1^s} - \fr 2/{2^s} \right)
			+ \left( \fr 1/{4^s} - \fr 2/{8^s} \right)
			+ \left( \fr 1/{9^s} - \fr 2/{18^s} \right)
			+ \left( \fr 1/{16^s} - \fr 2/{32^s} \right) + \cdots.
\\
\noalign{\text{After rearranging terms by absolute convergence and using (\ref{thmbeta}), this becomes }}
	&= \fr 1/{1^s} - \fr 2/{2^s} + \fr 1/{4^s} - \fr 2/{8^s} + \fr 1/{9^s} + \fr 1/{16^s} 
		- \fr 2/{18^s} + \fr 1/{25^s} - \fr 2/{32^s} + \fr 1/{36^s} + \cdots,
\\
		&= \suny \, \fr {\beta(n)}/{n^s},
\end{align*}
and finally, inserting the definition (\ref{defbeta}) of $\beta(n)$, \cite[eq. 3.4, 3.5]{Lee:2013}
\begin{equation} \label{lhs}
	(1-2^{1-s}) \zeta(2s) = \suny \, \sumn 
				\fr { (-1)^{\Omega(m)} (-1)^{l+1} }/{m^s l^s}, \qquad(l=\fr n/m, \Re s >\fr 1/2).
\end{equation}

\bigskip

For a different proof of the last equation, consider the classic Liouville arithmetic function 
$\lambda(n)=(-1)^{\Omega(n)}$ \cite[p. 184]{Sierpinski:1988} whose generating function
is \cite[p. 618]{Landau:1909}
\begin{equation} \label{lambda}
	 \suny \, \fr {\lambda(n)}/{n^s} = \fr {\zeta(2s)}/{\zeta(s)}  \qquad(s>1). 
\end{equation}
Landau uses Euler's product for $\zeta(s)$ with $s>1$ in his proof and notes that the exact domain of
convergence of the $\lambda$-series is one of the most important unsolved problems of number theory.

A Dirichlet product \cite[p. 670]{Landau:1909} then gives the right-hand side of (\ref{lhs})
and the uniqueness of the coefficients of a Dirichlet series yields immediately the expression 
(\ref{thmbeta}) for the arithmetic function $\beta(n)$~:
\begin{align*}
  (1-2^{1-s})\zeta(2s)&=\eta(s)\fr {\zeta(2s)}/{\zeta(s)}
	= \suny\,\fr {(-1)^{n+1}}/{n^s}\,\suny\,\fr {\lambda(n)}/{n^s}
\\
	&=\suny \, \sumn \fr { \lambda(m) (-1)^{1+\tfr n/m}}/{n^s}
\\
	&= \suny \, \fr {\beta(n)}/{n^s},	\qquad(s>1).
\end{align*}
Since both $\zeta(2s)$ and the $\beta$-series (a sum over $n^2$) are analytic for $\Re s >\fr 1/2$,
the previous equality can be extended from $s>1$ to $\Re s>\fr 1/2$, yielding a complete proof of
the equation (\ref{lhs}).

\section{Series of zero series}
If $s$ is a zero of $\eta(s)$ such that $\Re s>0$, then for any natural number $m$  \cite[p. 10]{Lee:2013}
$$
   \fr {\lambda(m)}/{m^s}\,\suly\,\fr {(-1)^{l+1}}/{l^s} = \fr {\lambda(m)}/{m^s}\,(0)=0,
 $$
which obviously yields \cite[p. 11]{Lee:2013}
\begin{equation} \label{rhs}
	\sumy\, \suly\,\fr {\lambda(m)(-1)^{l+1}}/{m^s l^s} = 0, \qquad (\Re s > 0, \eta(s)=0).
\end{equation}

\section{Comparison of double series}
In order to prove (\ref{Leeclaim}), Mr. Lee states in \cite[theorem 3.10]{Lee:2013} that the double
sums in (\ref{lhs}) and (\ref{rhs}) are equal because the product sets of integers $(n,l)$ in (\ref{lhs})
and $(m,l)$ in (\ref{rhs}) are identical. In fact, what needs to be justified is that the
order of summation in the double series can be changed~:
\begin{align}
\label{Snm}
(1-2^{1-s}) \zeta(2s) = & 
	\suny \, \sumn \fr { \lambda(m)(-1)^{l+1}}/{m^s l^s}, \qquad(l=\fr n/m)
\\ \label{Smn}
=& \sumy\, \sum_{n=m, m\mid n}^\infty \fr {\lambda(m)(-1)^{l+1}}/{m^s l^s}, \qquad(l=\fr n/m)
\\ \nonumber
=& \sumy\, \suly \fr {\lambda(m)(-1)^{l+1}}/{m^s l^s} \,=\, 0, \quad (\fr 1/2<\Re s <1, \eta(s)=0).
\end{align}
Some justification is needed here for the inversion between (\ref{Snm}) and (\ref{Smn}) since
according to Riemann's rearrangement theorem for single series \cite[p. 318]{Knopp:1954}, 
\begin{quote}
``If a series converges, but not absolutely, its sum can be made to have
any arbitrary value by a suitable derangement of the series; it can also be made divergent or
oscillatory."  \cite[p. 74]{Bromwich:1955}
\end{quote}
Although the original and the reordered series have exactly the same terms, the second 
can diverge or converge conditionally to a different sum.

A double series can be summed by columns (\ref{lhs},\ref{Snm}), 
by rows (\ref{rhs},\ref{Smn}), or by expanding rectangles in 
the sense of Pringsheim \cite[\S 2.5]{WhittakerWatson:1927}. 
For positive terms, the three sums are equal if any one of them converges.
\begin{quote}
``When the terms of the double series are positive, its convergence implies the convergence
of all the rows and columns, and its sum is equal to the sum of the two repeated series."
\cite[p. 84]{Bromwich:1955}

``The terms being always positive, if either repeated series is convergent, so is the other
and also the double series; and the three sums are the same." \cite[p. 84]{Bromwich:1955}
\end{quote}
These properties can be extended easily to absolutely convergent double series.
\begin{quote}
``It is clear that an array whose elements are indiscriminately positive and negative, if it converges
absolutely, may be treated as if it were a convergent array of positive terms."
\cite[p. 356]{Goursat:1904}
\end{quote}
But the situation for general terms is not so simple. 
A striking example from Ces\`aro makes this clear \cite[p. 89]{Bromwich:1955}~:
if $a_{m,n}=(-1)^{n+1}b_n(1-b_n)^{m-1}$ where $b_n=\tfr 1/{2^{\lfloor{n/2}\rfloor+1}}$,
then the sum of row $m$ is $\tfr 1/{2^m}$, converging absolutely, and so the 
sum by columns is 1. But the sum of column $n$ is $(-1)^{n+1}$, so the sum
by rows is oscillating. In a similar example due to Arndt \cite[p. 356]{Goursat:1904},
the sum by columns is $\tfr 1/2$ while the sum by rows is $-\tfr 1/2$.
Thus
\begin{quote}
``the sum of a non-absolutely convergent double series may have different values according
to the mode of summation"  \cite[p. 89]{Bromwich:1955}.
\end{quote}
\begin{quote}
``A double series should not be used in computations unless it is
 absolutely convergent. \cite[p. 357]{Goursat:1904}.
\end{quote}

\bigskip

However, Pringsheim has proven \cite[p. 28]{WhittakerWatson:1927} the following result~:
\begin{quote}
``If the rows and columns converge, and if the double series is convergent, then the
repeated sums are equal" \cite[p. 81]{Bromwich:1955}.
\end{quote}

For the double series in (\ref{Smn}), column $n$ has a finite number of terms whose sum is $\tfr \beta(n)/n^s$,
while the sum of each row $m$ converges to $0$ if $\eta(s)=0$. But it is not proven in 
\cite{Lee:2013} that the sums by expanding rectangles converge, and
Pringsheim's theorem cannot be used to show that the sum of the double series in (\ref{Snm}) is zero.

\section{Refutation of the first version of Mr Lee's proof}
The simple proof of the ``Riemann Hypothesis" proposed in \cite{Lee:2013}, although interesting 
and original, is clearly incomplete~: a crucial theorem presents conditionally convergent 
infinite series as sums over sets, without specifying the order of summation, and without 
providing any justification for disregarding this order.

\section{Uniform convergence and double sequences}
\subsection{}
After being  made aware of this gap in his proof, the author of \cite{Lee:2012} and \cite{Lee:2013}
suspended \cite{Lee:2012} and proposed in \cite{Lee:2013a} a justification based on the Moore 
theorem for the inversion of two limits, one of which is uniform \cite[p. 28]{DundordSchwartz:1957}.
The second version of Mr Lee's proof was made public one year later in \cite{Lee:2014}, and relies
on two theorems (2.13, 2.15) from a very readable elementary report about double sequences made
available over the Internet in 2005 by Dr Eissa Habil \cite{Habil:2005}.

\subsection{}
Below, the expression ``exist-U" means that the convergence is uniform with respect 
to the free variable in $\mathbb{N}$, and we will use the following abbreviations~:
$$
\begin{array}{l|ll}
	\text{Symbol} & \text{Long notation} & \text{Description}
\\
\hline
\\
	f(\infty,n) & \lim_{m\to\infty} f(m,n) & \text{first partial limit}, n\in\mathbb{N}
\\
	f(m,\infty) & \lim_{n\to\infty} f(m,n) & \text{second partial limit}, m\in\mathbb{N}
\\
	f(\infty,\infty) & \lim_{m,n\to\infty} f(m,n) & \text{limit of double sequence}
\\
	f(\infty_1,\infty_2) & \lim_{m\to\infty}\lim_{n\to\infty} f(m,n) & \text{first iterated limit}
\\
	f(\infty_2,\infty_1) & \lim_{n\to\infty}\lim_{m\to\infty} f(m,n) & \text{second iterated limit}
\end{array}
 $$

\subsection{}
Following \cite[Definition 2.1]{Habil:2005}, a double sequence $f(m,n)$ of complex numbers
converges to zero if and only if
$$
	(\forall \epsilon>0) (\exists N \in \mathbb{N}) \text{ such that }
		m>N \And n>N \Rightarrow |f(m,n)|<\epsilon.
 $$

\subsection{}\label{thm215}
Two conditions sufficient for such convergence to zero are specified in 
\cite[Theorem 2.15]{Habil:2005}~:
$$
	f(m,\infty) \text{ exists-U } \And f(\infty_1,\infty_2)=0 \implies f(\infty,\infty)=0.
 $$
Proof: In the inequality
$$
\left|f(m,n)\right| \le  \left|f(m,n) - f(m,\infty)\right| + \left|f(m,\infty)\right|,
 $$
the last term is small if $m$ is large by hypothesis, while the middle difference
is small if $n$ is large, independently of the previously chosen large $m$ by uniformity.

\subsection{}\label{thm211c}
Conversely \cite[Theorem 2.11, corrected]{Habil:2005}~:
$$
	f(m,\infty) \text{ exists-U } \And f(\infty,\infty)=0 \implies f(\infty_1,\infty_2)=0.
 $$
Proof:  In the inequality
\begin{equation} \label{ineqfmn}
\left|f(m,\infty)\right| \le \left|f(m,\infty) - f(m,n)\right| + \left|f(m,n)\right|,
\end{equation}
the last term is small if $m$ and $n$ are large by hypothesis, while the middle difference
is small if $n$ is large, independently of the previously chosen large $m$ by uniformity.

\subsection{}\label{thm213c}
As a corollary \cite[Theorem 2.13, corrected]{Habil:2005}~:
$$
 f(\infty_1,\infty_2)=0 \And f(m,\infty) \text{ exists-U } \And f(\infty,n) \text{ exists-U } 
		\implies  f(\infty_2,\infty_1)=0.
 $$
Moore's theorem \cite[p. 28]{DundordSchwartz:1957} is stronger than this, requiring only {\em one} uniform limit~:
$$
 f(m,\infty) \text{ exists } \And f(\infty,n) \text{ exists-U } 
		\implies  f(\infty,\infty)=f(\infty_1,\infty_2) = f(\infty_2,\infty_1).
 $$

\subsection{}
Unfortunately, the important uniformity condition has been omitted in \cite[Theorems 2.11, 2.12, 2.13]{Habil:2005},
for example in theorem (2.11)~:
$$
	f(m,\infty) \text{ \em exists } \And f(\infty,\infty)=0 \implies f(\infty_1,\infty_2)=0.
 $$
The proposed proof relies on the inequality (\ref{ineqfmn}) and on the three bounds:
\\
1. ``Given $\epsilon>0$, there exists $N_1$ such that $|f(m,n)|<\tfr\epsilon/2$ if $m,n>N_1$";
\\
2. ``There exists $N_2$ such that $|f(m,\infty)-f(m,n)|<\tfr\epsilon/2$ if $n>N_2$";
\\
3. ``Choose $n>\max(N_1,N_2)$. Then $\forall m>N_1, |f(m,\infty)| \le \epsilon$".
\\
In general, the integer $N_2$ could however depend on $m$ in the second bound, so that it could be impossible
to choose $n>N_2, \forall m>N_1$ for the third bound! Fortunately, the hypothesis can be modified
easily by replacing ``{\em exists}" with ``{\em exists uniformly}" in order to make the proof correct,
as was done in \S\ref{thm211c}. This stronger condition already appears in the statement
of \cite[Theorem 2.15]{Habil:2005} in \S\ref{thm215}.

\section{Uniform convergence and double series}
\subsection{}
The convergence to zero of a double series of complex numbers $a_{m,n}$ is defined 
as the convergence to zero of the double sequence of partial sums
\cite[Definition 7.1]{Habil:2005}
$$
	S(M,N)=\sum_{m=1}^M\sum_{n=1}^N a_{m,n}, \qquad(M,N\in\mathbb{N}).
 $$

\subsection{}
Following \cite[Definition 3.1]{Lee:2014}, let
$$
	a_{m,n} =	\begin{cases}
				\tfr {\lambda(m)(-1)^{1+\tfr n/m}}/{n^s},	&\text{if }m\text{ divides } n
\\
				0							&\text{otherwise}
			\end{cases}
 $$
where $\tfr 1/2<\Re s<1$ and $\eta(s)=0$, if possible.
The iterated series in (\ref{Snm}) then corresponds to $S(\infty_2,\infty_1)$, 
and the iterated series in (\ref{Smn}) to $S(\infty_1,\infty_2)$.

\subsection{}
Both partial limits certainly exist in this case, from (\ref{lhs}) and (\ref{rhs})~:
$$
S(\infty,N)=\sum_{n=1}^N \fr {\beta(n)}/{n^s}, \forall N\in\mathbb{N};\quad
S(M,\infty)=\sum_{m=1}^M \fr {\lambda(m)}/{m^s}\eta(s)=0, \forall M\in\mathbb{N}.
 $$
In \cite[Lemma 3.5]{Lee:2014}, only the differences 
$S(\infty,n)-S(\infty,n-1)=\tfr {\beta(n)}/{n^s}$
and $S(m,\infty)-S(m-1,\infty)=0$ are verified, but this is equivalent by 
telescoping finite sums. Obviously, 
$S(\infty_1,\infty_2)=\lim_{M\to\infty}S(M,\infty)=0$ also.

\subsection{}
From the stronger version of corollary \ref{thm213c} derived by Moore's theorem,
$$
 S(\infty_1,\infty_2)=0 \And S(M,\infty) \text{ exists } \And S(\infty,N) \text{ exists-U } 
		\implies  S(\infty_2,\infty_1)=0.
 $$
In order to prove that $(1-2^{1-s})\zeta(2s)=S(\infty_2,\infty_1)=0$,
we thus only need the uniformity with respect to $N\in\mathbb{N}$ 
in the partial limit
$$
  S(\infty,N) = 
	\lim_{M\to\infty} \sum_{m=1}^M \sum_{n=m,m\mid n}^N 
		\fr {(-1)^{\Omega(m)}(-1)^{1+\tfr n/m}}/{n^s}.
 $$
Equivalently \cite[p. 123]{Bromwich:1955}, we must show that
$$
	(\forall \epsilon>0) (\exists M \in \mathbb{N}) \text{ such that }
		\left|\sum_{m=M}^{M+p} \sum_{n=m}^N a_{m,n}\right|<\epsilon,
		 \qquad(\forall N\in\mathbb{N}, \forall p\in\mathbb{N}).
 $$

\subsection{}
Unfortunately, in \cite[Lemma 3.5]{Lee:2014}, we only find the verification
of a different result, ``$\suny a_{m,n}$ converges to zero uniformly on $m$", 
which is equivalent to
$$
	(\forall \epsilon>0) (\exists N \in \mathbb{N}) \text{ such that }
		 \left|\sum_{n=N}^{N+p} a_{m,n}\right|<\epsilon,
		 \qquad(\forall m\in\mathbb{N}, \forall p\in\mathbb{N}).
 $$

\section{Refutation of the second version of Mr Lee's proof}

The simple proof of the ``Riemann Hypothesis" proposed in \cite{Lee:2014} is incomplete~: it refers 
to an unproven, probably misstated, theorem from another report that has not been published. 
Even if this crucial theorem was accepted as is, or more appropriately restated to 
make its current proof correct, only two of three needed conditions are verified in \cite{Lee:2014}.
The gap found in the first version of the proof \cite{Lee:2013} has not been filled in the
second version \cite{Lee:2014}.

\bibliographystyle{plain}

\end{document}